\documentclass[12pt]{article}
\usepackage{graphicx}

\font\tendb=msbm10 at 12pt \font\sevendb=msbm10 at 9pt
\font\fivedb=msbm10 at 7pt
\newfam\dbfam
\textfont\dbfam=\tendb \scriptfont\dbfam=\sevendb
\scriptscriptfont\dbfam=\fivedb
\def\db{\fam\dbfam\tendb}
\font\eufm=eufm10\font\eufms=eufm10\font\eufmss=eufm10\newfam\eufam
\textfont\eufam=\eufm\scriptfont\eufam=\eufms\scriptscriptfont\eufam=\eufmss

\font\tendbb=msbm10 at 12pt \font\sevendbb=msbm7 at 9pt
\font\fivedbb=msbm5 at 6pt
\newfam\dbbfam
\textfont\dbbfam=\tendbb \scriptfont\dbbfam=\sevendbb
\scriptscriptfont\dbbfam=\fivedbb

 \def \Z {{\db Z}}
 \def \R {\hbox{\db R}}
 
 \def \C {\hbox{\db C}}
 
 \def \S {S^{3}}

\font\tenMmm=eusm10 at 12pt
\newfam\Mmmfam
\textfont\Mmmfam=\tenMmm

\def\illu #1 by #2 (#3){
  \vbox to #2{
    \hrule width #1 height 0pt depth 0pt
    \vfill
    \special{illustration #3} % this is the low-level interface
    }
  }

\begin{document}

\null \vspace{3cm}

\begin{center}
{\large {\bf The Casson-Walker-Lescop  invariant of periodic three-manifolds}}\\
 Nafaa Chbili\footnote{ Supported by a  fellowship from the Japan Society for the
Promotion of Science (JSPS), and by a Grant-in-Aid for JSPS fellows 03020.}\\
 \begin{footnotesize} Department of Mathematics\\
 Tokyo Institute of Technology\\
 Oh-okayama Meguro Tokyo 152-8551 JAPAN\\
 E-mail: chbili@math.titech.ac.jp
 \end{footnotesize}
\end{center}
\begin{footnotesize}
               {\bf Abstract.} Let $p$ be an odd prime and $G$ the finite cyclic group of order $p$. We use the
                Casson-Walker-Lescop  invariant
               to find a necessary condition for a three-manifold
               to have an action of $G$ with a circle as the set
               of fixed points.

               \end{footnotesize}

\begin{center}
{\sc I- Introduction}\\
\end{center}

 In 1985, A. Casson introduced an integer valued invariant of oriented integral homology spheres which  reduces modulo
  two to the Rokhlin invariant. In some sens, the Casson
invariant counts the number of conjugacy classes of
representations of the fundamental group  into $SU(2)$. K. Walker
\cite{Wa} extended the Casson invariant to rational homology
spheres and established a method of defining this invariant, that
relies only on the surgery presentation of the considered
manifold. By a purely combinatorial approach, C. Lescop \cite{Les}
extended the Casson-Walker invariant to any oriented 3-manifold.
Namely, given a 3-manifold $M$ and $L$ a surgery presentation of
$M$, Lescop introduced an explicit formula for the Casson-Walker
invariant defined in terms of the Alexander polynomials, linking
numbers and framing coefficients of the sub-links of L. Lescop's
formula will be our
main tool in this paper.\\
Let $M$ be an oriented closed 3-manifold and $G$ a finite group.
An intriguing question in low dimensional topology is to know
whether the group $G$ acts non trivially on $M$. Many different
approaches were applied to study this question. In the last
decade, the  quantum invariants were used in order to find
necessary conditions for a three-manifold to have such an action,
see \cite{BGP}, \cite{Ch1}, \cite{Ch2}, \cite{CL} and  \cite{Gi}.
The main motivation of this approach is the fact that the quantum
invariants of links were used in a successful way to study
symmetries of links and knots in the three-sphere.\\
 A link $L$ in
$\S$ is said to be $p$-periodic if there exists an
orientation-preserving  periodic transformation $h$ of order $p$
such that $fix(h)\cong S^1, h(L)=L$ and   $fix(h)\cap
L=\emptyset$. By the positive solution of the Smith conjecture
\cite{BM}, we may assume that $h$ is a rotation by a $2\pi/p$
angle around the $z$-axis. An oriented closed  3-manifold is said
to be $p$-periodic if the finite group $\Z/p\Z$ acts semi-freely
on
$M$, with a circle as the set of fixed points.\\
A $p$-periodic link is said to be strongly  $p-$periodic  if  each
component of the quotient link $\overline L$ has linking number
zero modulo $p$ with the axis of the rotation. A strongly periodic
link $L$ is said to be orbitally separated if the quotient link is
algebraically split. In \cite{PS}, Przytycki and Sokolov proved
that a three-manifold is $p$-periodic if and only if $M$ is
obtained from the three-sphere by surgery along a strongly $p$-periodic link.\\
%It is
%a natural idea to use the periodic surgery presentation introduced
%in \cite{PS} to study the quantum invariants of periodic
%three-manifolds. This was the approach in
There are several studies about the Casson invariant of cyclic
branched covers. In particular, Davidow \cite{Da}, Hoste
\cite{Ho1} and Ishibe \cite{Is}, studied independently the case of
Whitehead doubles of knots. Mullins \cite{Mu}, considered the case
of 2-fold branched covers along links and provided a nice formula
relating the Casson-Walker invariant to the 2-signature of the
link and some derivatives of the Jones polynomial. In \cite{GK},
 Garoufalidis and Kricker computed the LMO invariant of cyclic
 branched covers and gave a formula for the Casson-Walker
 invariant under the assumption that the branched cover is a
 rational homology three-sphere.\\
 The purpose of the present paper
is to study the Casson-Walker-Lescop  invariant of periodic
three-manifolds. The only restriction we need here is that   the
quotient manifold is an integral homology three-sphere. In other
words, we deal with manifolds which arise as cyclic branched
covers, along knots, of  integral homology spheres.\\
We first study  the second coefficient of the Conway polynomial of
strongly periodic links. Then we use the  surgery formula
introduced by Lescop to find a necessary condition for a
three-manifold to be periodic.\\

{\bf Notations:} Let $M$ be an oriented closed three-manifold. We
denote by $\lambda(M)$ the Casson-Walker-Lescop invariant of $M$.
It is worth mentioning  that $\lambda$ takes values in
$\frac{1}{12}\Z$. Let $L=l_1\cup\dots\cup l_n$ be a $n-$component
framed link. For $1\leq i \leq n$,  let $l_{ii}$ be the framing
coefficient  of the component $l_i$ and define $l_{ij}$ to be the
linking number of $l_{i}$ and $l_{j}$ if $i \neq j$. Let
$E(L)=(l_{ij})$ be  the symmetric linking matrix of $L$. By
$Sig(E(L))$ we denote the signature of $E(L)$. Finally, $|H_1(M)|$
denotes the number of elements of  the first homology group
$H_1(M)$ if finite and zero otherwise.
Here are the main results that we will prove in this paper.\\

{\bf Theorem 1.} {\sl Let $p\geq 3$ be a prime and $M$ a
$p$-periodic $3$-manifold such that the quotient $\overline M$ is
an integral homology sphere. Then there exists an orbitally
separated strongly $p$-periodic link such that $M$ is obtained
from $\S$ by surgery along $L$ and :
$$24\lambda(M)\equiv 3|H_1(M)|Sig(E(L)) \mbox{ modulo } p,
$$
}

{\bf Corollary 1.} {\sl Let $p\geq 3$ be a prime and $M$ a
$p-$periodic manifold. If the first Betti number of $M$ is not
zero, then
$24\lambda(M)\equiv 0 \mbox{ modulo } p.$}\\

Let $M$ be an integral  homology sphere. We know that $M$ can be
obtained from $\S$ by  surgery along an algebraically split link.
On the other hand, we know that if $M$ is $p$-periodic then $M$
has a strongly $p$-periodic orbitally separated surgery
presentation. One may ask if a periodic homology sphere can be
obtained from $\S$ by surgery along an
algebraically split strongly periodic link.\\

{\bf Corollary 2.} {\sl Let $p\geq 3$ be a prime and $M$ a
$p$-periodic 3-manifold. If $M$ has a strongly $p$-periodic
algebraically split  surgery presentation
then $24\lambda(M)\equiv 0 \mbox{ modulo } p.$ }\\

{\bf Example.} Let $T^3=S^1 \times S^1 \times S^1$ be  the
three-dimensional torus. The Casson  invariant of $T^3$ is equal
to one. By applying Corollary 1, we can see easily that $T^3$ is
not $p$-periodic for all prime
$p>3$. Same conclusion for the manifold $S^1 \times S^2$. \\

If $L$ is a link with $n$ components, let $\nabla_L$ denote the
Conway polynomial (see paragraph 3) of $L$, we know that  $
\nabla_{L}(z)=z^{n-1}(a_0 + a_1 z^2+ \dots + a_m z^{2m})$, where
coefficients $a_i$ are integers. The proof of Theorem 1 relies
heavily on the following
theorem, which will be proved in paragraph 5.\\
{\bf Theorem 2.} {\sl Let $p \geq 3$ a prime and $L$ an orbitally
separated strongly $p$-periodic link. Then the coefficient $a_1$
is null modulo
$p$.}\\

By considering the trefoil knot which is 3-periodic, we can see
easily that the condition  given by  the Theorem 2 is not
satisfied by periodic links in general. Indeed,  a simple
computation shows that for the trefoil knot  the coefficient
$a_1=1$, hence it is not null modulo
$3$.\\
This paper is outlined as follows. In section 2, we introduce
strongly periodic links and some of their properties. In section
3, we review  the surgery presentation of periodic manifolds.
Section 4 is to review  some properties of the Conway polynomial
needed in the sequel. The proof of Theorem 2 is given in section
5. In section 6, we briefly recall  Lescop's surgery formula.
Ultimately, in section 7, we give the proof of Theorem 1.\\

%%%%%%%%%%%%%%%%%%%%%%%%%%%%%%%%%%%%%%%%%%%%%%%
%%%%%%%%%%%%%%%%%%%%%%%
\begin{center}
{\sc II- Strongly periodic links.} \end{center}
%%%%%%%%%%%%%%%%%%%%%%%%%%%%%%%%%%%%%%%%%%%%%%%%%%%%%%%%%%%%%%%%%%%%%%
In this paragraph we recall some definitions.\\
{\bf Definition 2.1.} {\sl { Let $p\geq 2$ be an integer. A link
$L$ of $\S$ is said to be $p$-periodic if and only if there exists
an orientation-preserving  auto-diffeomorphism  of  $\S$
such that:\\
1- Fix($h$) is homeomorphic to the  circle $S^{1}$,\\
2- the link  L is disjoint from  Fix($h$),\\
3- $h$ is of order $p$,\\
4- $h(L)=L$.\\
If $L$ is  periodic we will denote the quotient link by $\overline
L$.\\}}

The standard example of such a diffeomorphism is given as follows.
Let us consider  $\S$ as the  sub-manifold  of $\C ^2$ defined by
$\S=\{(z_{1},z_{2}) \in \C ^{2}; |z_{1}|^{2}+|z_{2}|^{2}=1\}$ and
$\varphi{_p}$ the following diffeomorphism:
$$
\begin{array}{cccl}
\varphi_{p}:& S^{3} & \longrightarrow & S^{3} \\
   & (z_{1},z_{2}) & \longmapsto & (e^{\frac{2i\pi}{p}}z_{1},z_{2}).
\end{array} $$
The set of fixed points for  $\varphi_{p}$ is the circle $\Delta =
\{(0,z_{2}) \in \C ^{2}; |z_{2}|^{2}=1\}$. If we identify $\S$
with $\R ^{3}\cup {\infty}$, $\Delta$ may be seen as the standard
$z$-axis. By the positive solution of the Smith conjecture
\cite{BM}, a periodic diffeomorphism such as in the previous
definition is conjugate to the standard rotation $\varphi_{p}$.
Consequently, if $L$ is a
$p$-periodic link then there exists a tangle $T$ such that $L$ is the closure of $T^p$.\\
%\vspace{3cm}
%\begin{picture}(0,0)
%\put(210,0){\framebox(40,20){\footnotesize {T}}}
%\put(210,-20){\framebox(40,20){\footnotesize {T}}}
%\put(230,-25){.} \put(230,-30){.} \put(230,-35){.}
%\put(210,-60){\framebox(40,20){\footnotesize {T}}}
%\end{picture}
Recall here that if the quotient link $\overline L$ is a knot then
the link $L$ may have more than one component. In general, the
number of components of $L$ depends on the linking numbers of
components of $\overline L$ with the axis
of the rotation.\\
Recently, Przytycki and Sokolov \cite{PS} introduced the notion of strongly periodic links as follows.\\

{\bf Definition 2.2.} {\sl { Let $p\geq 2$ be an integer. A
$p$-periodic link $L$ is said to be strongly $p$-periodic if and
only if one of the following conditions holds:\\
i) The linking number of each component of $\overline L$ with the
axis $\Delta$ is congruent
 to zero  modulo $p$.\\
 ii) The group $\Z/{p\Z}$ acts freely on the set of components of
 $L$.\\
 iii) The number of components of $L$ is $p$ times greater than
 the number of components $
\overline L$. }} \\
 {\bf Remark 2.1.} According to condition iii)  in the previous definition, a $p$-strongly
 periodic link $L$ has $p\alpha$ components, where $\alpha$ is the
 number of components of the quotient link $\overline L$. These
 $p\alpha$ components are divided into $\alpha$ orbits with
 respect to the free action of  $\Z/{p\Z}$ (condition i)).\\
 Assume that  $\overline L=l_1\cup \dots \cup\l_{\alpha}$, there is a natural
cyclic order on the components of $L$. Namely,
$$L=l_1^1 \cup\dots \cup l_1^p \cup l_2^1 \cup\dots \cup l_2^p \cup \dots \cup l_{\alpha}^1 \cup\dots
\cup l_{\alpha}^p$$ where $ \varphi_p(l_i^t)=l_{i}^{t+1} \;
\forall 1\leq t\leq p-1$ and $\varphi_p(l_i^p)=l_{i}^{1}$, for all
$1\leq i \leq \alpha$.

{\bf Example.} The link in the following figure is a strongly
5-periodic
link.\\
 \begin{center}
\includegraphics[width=5cm,height=3cm]{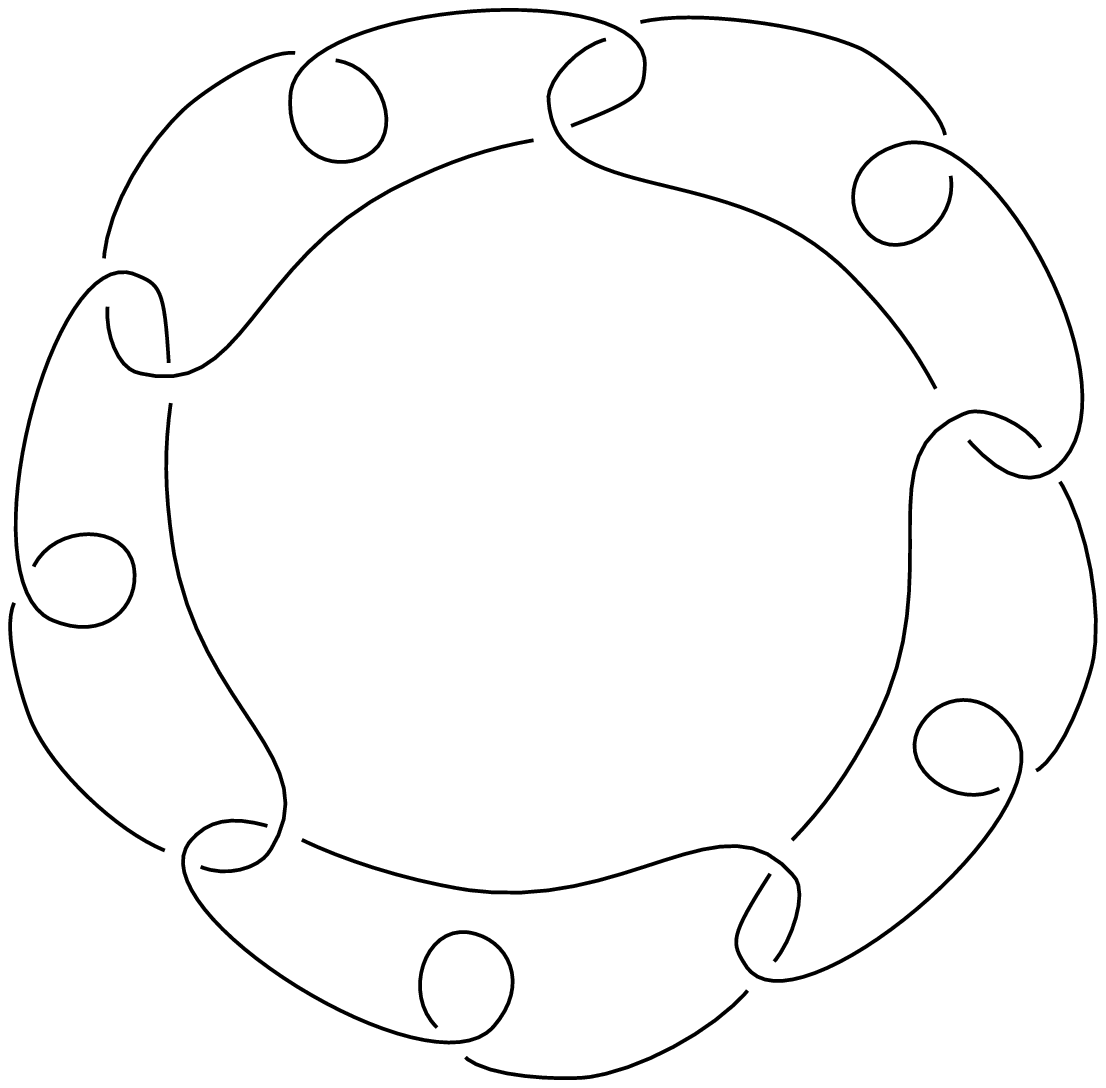}
\end{center}
\begin{center}{\sc  Figure 1. }\end{center}
 {\bf Definition 2.3.}  {\sl Let $L$ be a link in the three-sphere.
 $L$ is said to be algebraically split if and only if the linking number of
 any two components of $L$ is null.}\\
 {\bf Definition 2.4.} {\sl { Let $p\geq 2$ be an integer. A strongly $p$-periodic link is said to be
 orbitally separated if and only if the quotient link is
 algebraically split.}}\\
 {\bf Remark 2.2.} Throughout the rest of this paper we use the term  OS link to refer to an orbitally separated link.
  If $L$ is a  strongly $p$-periodic OS link, then for all $s$ and $i\neq j$ we have
  $\displaystyle\sum_{t=1}^p lk(l_i^s,l_j^t)=0$.

%%%%%%%%%%%%%%%%%%%%%%%%%%%%%%%%%%%%%%%%%%%%%%%%%%%%%%%%%%%%%%%%%%%%%%%%%%%%%%%%%%%%%%%%%%%%%%%%%%%%%%%%%%%%%%%%%%%%%%%%%%%%%%%%%%%%%%%%%%%%%%%%%%%%%%%%%%%%%%%%%%%%%%%%%%%%%%%%%%%%%%%%%%%%%%%%%%%%%%%%%%%%%%%%%%%%%%%%%%
\begin{center}{\sc III- Periodic 3-manifolds }\end{center}

Let $M$ be a 3-manifold. It is well known from Lickorish's
 work \cite{Li1} that $M$ may be obtained from the three-sphere  $\S$ by
 surgery along a framed link $L$ in $\S$. Such a link is called a surgery
  presentation of $M$.
    This section deals  with surgery presentations of periodic three-manifolds. \\
 {\bf Definition 3.1.} {\sl { Let $p\geq 2$ be an integer.  A
3-manifold $M$ is said to be $p$-periodic
if and only if there exists an orientation-preserving  auto-diffeomorphism of $M$, such that:\\
1- Fix($h$) is homeomorphic to the  circle $S^{1}$,\\
2- $h$ is of order $p$.\\ }}

The three-sphere $\S$ is $p$-periodic for all  $p\geq 2$.  The
acting diffeomorphism is nothing but the
 rotation $\varphi_{p}$ defined in section 2. Another  concrete example is the Brieskorn manifold
 $\Sigma(p,q,r)$,
 which has periods $p$, $q$ and $r$. One can easily give explicit  formulas for   the corresponding diffeomorphisms by
  considering $\Sigma(p,q,r)$ as the intersection of the complex surface $\{(z_{1},z_{2},z_{3} )
 \in \C ^{3}; z_{1}^{p}+z_{2}^{q}+z_{3}^{r}=0\}$ with an appropriate sphere of dimension 5
  (see \cite{Mil}).\\

 {\bf Theorem 3.2 (\cite{PS}).} {\sl Let $p$ be a prime. A 3-manifold  $M$ is  $p$-periodic if and only if
  $M$ is obtained from
$\S$ by surgery along a strongly $p-$periodic link. }\\
{\bf Remark.} If the quotient manifold is a homology sphere then
we may suppose that  $L$ is orbitally separated.\\

%%%%%%%%%%%%%%%%%%%%%%%%%%%%%%%%%%%%%%%%%%%%%%%%%%%%%%%%%%%%%%%%%%%%%%%%%%
%%%%%%%%%%%%%%%%%%%%%%%%%%%%%%%%%%%%%%%%%%%%%%%%%%%%%%%%%%%%%%%%%%%%%%%%%%

\begin{center}
{\sc IV- The Conway polynomial}\\
\end{center}

The Alexander polynomial of periodic knots has been subject of
 an extensive literature. In 1971, Murasugi \cite{Mu} proved a congruence
 relationship between the Alexander invariant of a periodic knot
 $K$ and the Alexander invariant of the quotient knot. Later,
 different approaches   were used to give alternative  proofs to Murasugi's
 result \cite{Bu}, \cite{DL},  \cite{Hi}. In \cite{Sa} and \cite{Mi}, similar
 results were proved for the multi-variable Alexander
 polynomial. Murasugi's result was applied to rule out the
 possibility of being periodic
 for many knots.\\
 In this paper we consider the Conway version of the Alexander
 polynomial. The Conway polynomial $\nabla$ is an invariant of ambient isotopy
of oriented links which can be defined uniquely by the following:
$$\begin{array}{ll}
 {\bf (i)}&\nabla_ {\bigcirc}(z)=1\\
{\bf (ii)}&\nabla_{L_{+}}(z)-\nabla_{L_{-}}(z)
=z\nabla_{L_{0}}(z),
\end{array}$$
where  $\bigcirc $ is the trivial knot, $L_{+}$,  $L_{-}$ and
$L_{0}$ are three oriented links which are identical except
 near one crossing where they look
like in the following figure:
\begin{center}
\includegraphics[width=8cm,height=2cm]{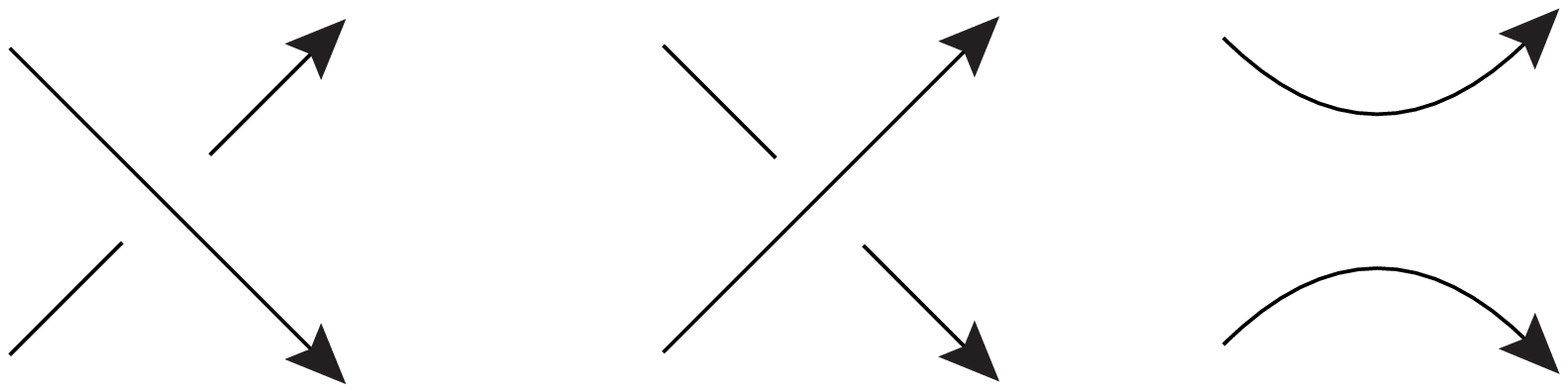}
\end{center}
$\hspace{140pt} L_{+}\hspace{80pt} L_{-}\hspace{70pt} L_{0}$
\begin{center} {\sc Figure 2}
\end{center}

It is well known that if $L$ is a link with $n$
 components then the Conway polynomial of $L$ is of the form $
\nabla_{L}(z)=z^{n-1}(a_0 + a_1 z^2+ \dots + a_m z^{2m})$, where
coefficients $a_i$ are integers. The first coefficient depends
only on the linking matrix of the link. The second coefficient
$a_{1}$ is the first non trivial Vassiliev invariant for knots. It
is related to the Casson-Walker invariant
by the global surgery formula introduced by Lescop \cite{Les}.\\
At the beginning, let us  recall some properties of the Conway
polynomial needed in the
sequel.\\
Let $L=l_1\cup \dots \cup l_n$ be an oriented $n$ component link
in the three-sphere. Let $l_{ij}=\mbox{Lk}(l_i,l_j)$ if $i \neq j$
and define $l_{ii}=-\displaystyle\sum_{j=1,j \neq i}^{n} l_{ij}$.
Define the linking matrix  ${\cal L}$ as ${\cal L}=(l_{ij})$. The
matrix $\cal L$ is symmetric and  each row adds to zero. Hence,
every cofactor of this matrix is the same. We refer the reader to
\cite{Ha}, \cite{Hoso} and \cite{Ho}, for
details about  the linking matrix.\\
 {\bf Theorem 4.3 \cite{Ho}.} {\sl Let $L$ an oriented link with
$n$ components. Then $a_0={\cal L}_{ij}$, where  ${\cal L}_{ij}$
is
any cofactor of the the linking matrix ${\cal L}$.}\\

It follows from this theorem, that if $L$ is algebraically split
then the first coefficient of the Conway polynomial is zero.
Levine \cite{Le} proved that the coefficient of $z^i$ is null, for
all $i\leq 2n-3$, and
gave an explicit formula for the coefficient of $z^{2n-2}$.\\
 {\bf Proposition 4.4 \cite{Le}.} {\sl If  $L$ is algebraically split
with $n$ components, then $\nabla_{L}(z)$ is divisible by
$z^{2n-2}$.}\\

Let $\Sigma_1$, $\Sigma_2$ and  $\Sigma_3$ be  three oriented
surfaces imbedded in $\S$. For $1\leq i \leq 3$, let $n_i$ be a
positive normal vector to $\Sigma_i$. Assume that the three
surfaces intersect transversally in their interiors. Let
$A=\{a_1,\dots,a_s\}=\Sigma_1 \cap \Sigma_2 \cap \Sigma_3$, every
$a_i \in A$ may be occupied  with a sign $\epsilon(a_i)=\pm 1$,
depending on the orientation of the basis defined by the three
normal vectors. Let $\Sigma_1:\Sigma_2:
\Sigma_3=\sum_{i=1}^{s}\epsilon(a_i)$.\\
Let $L=l_1 \cup\dots \cup l_n$ be an algebraically split link. We
can choose Seifert surfaces $\Sigma_1,$ ... $,\Sigma_n$ such that
$l_i$ bounds  $\Sigma_i$  and $l_j$ does not intersect $\Sigma_i$
if $i \neq j$. We define  $\bar \mu _{ijk}$ to be
$\Sigma_1:\Sigma_2: \Sigma_3$ if $i$, $j$ and  $k$ are distinct
and to be zero otherwise.  Now we define $u$ to be the
$(n-1)\times (n-1)$-matrix defined by $u_{ij}=\sum_{k=1}^{n} \bar
\mu_{ijk} ,$ for every $1\leq i,j \leq n-1$. Levine \cite{Le},
proved the following.\\
{\bf Theorem 4.5 \cite{Le}.} {\sl If  $L= l_1 \cup\dots \cup l_n$
is an algebraically split link then: $ a_{n-1}=det(u). $ }\\

%%%%%%%%%%%%%%%%%%%%%%%%%%%%%%%%%%%%%%%%%%%%%%%%%%%%%%%%%%%%%%%%%%%%%%%%%%%%%%%%%%%%%%%%%%%%%%%%%
%%%%%%%%%%%%%%%%%%%%%%%%%%%%%%%%%%%%%%%%%%%%%%%%%%%%%%%%%%%%%%%%%%%%%%%%%%%%%%%%%%%%%%%%%%%%%%%%%
\begin{center}
{\sc V- Proof of Theorem 2}
\end{center}

This paragraph is devoted to the proof of Theorem 2. Let $L$ be a
strongly $p$-periodic link in the three-sphere. Let $\overline L$
be the factor link, so here we have $L=\pi^{-1}(\overline L)$,
where $\pi$ is the canonical surjection corresponding to the
action of the rotation on the three-sphere. Let $\overline L_+$,
$\overline L_-$ and $\overline L_0$ denote the three links which
are identical to $\overline L$ except near one crossing where they
are like in figure 2. Now, let $L_{p+}:=\pi^{-1}(\overline L_+)$,
$L_{p-}:=\pi^{-1}(\overline L_-)$ and
$L_{p0}:=\pi^{-1}(\overline L_0)$.\\

{\bf Lemma 5.1.} {\sl Let $p$ be a prime, then we have the
following congruence modulo $p$: }
$$
\nabla_{L_{p+}}(z)-\nabla_{L_{p-}}(z)\equiv z^p
\nabla_{L_{p0}}(z).
$$
{\bf Proof.} Przytycki \cite{Pr} proved a similar lemma for the
skein (HOMFLY) polynomial of periodic links. The same argument can
be used here to prove our lemma.\\

Notice that if $L_{p+}$ is a strongly $p$-periodic link then
$L_{p-}$ is also a strongly periodic link. However, two cases are
to be
considered for the link $L_{p0}$.\\
First case: if the considered crossing of the factor link is a
mixed crossing, (involving two different components of the factor
link say $K_1$ and $K_2$).  After the crossing change, the two
components connect and give birth to a new component $K_{12}$, as
the linking numbers of $K_1$ and $K_2$ with the rotation axis is
null modulo $p$ then the linking number of $K_{12}$ with the
rotation axis is also null modulo $p$ as it is the sum of two
zeros. Consequently, the link $L_{p0}$ is  strongly $p$-periodic.\\
Second case: If the considered crossing is a self crossing, in
this case the link $L_{p0}$ may be not  strongly $p$-periodic.\\
Let $L$ be  a strongly $p$-periodic link, then we know that $L$
has $p \alpha$ components, where $\alpha $ is a positive integer.
The proof of Theorem 2, will be done by induction on $\alpha$.
Recall here that $\alpha$ is the number of components of
the factor link $\overline L$.\\
It is easy to see that  each of $L_{p+}$ and $L_{p-}$ has
$p\alpha$ components. Let
$s$ be the number of components of  the link $L_0$. Let\\
$\nabla_{L_{p+}}(z)=z^{p\alpha-1}(a_0 + a_1 z^2+ \dots + a_m z^{2m}),$\\
$\nabla_{L_{p-}}(z)=z^{p\alpha-1}(b_0 + b_1 z^2+ \dots + b_k z^{2k}),$\\
$\nabla_{L_{p0}}(z)=z^{s-1}(c_0 + c_1 z^2+ \dots + c_r z^{2r}).$\\

{\bf Lemma 5.2.} {\sl Let $p\geq 3$ be a prime and  $L$ a strongly
$p$-periodic link. Assume that $L$ has $p$ components, then:} $a_1
(L)\equiv 0 \mbox{ modulo } p.$\\
 {\bf Proof.} Here  $\alpha=1$, then
each of $L_{p+}$ and $L_{p-}$ has $p$ components. By lemma 5.1, we
have the following congruence modulo $p$:
$$ z^{p-1}(a_0 + a_1 z^2+ \dots + a_m z^{2m})-z^{p-1}(b_0 + b_1 z^2+ \dots + b_k
z^{2k})\equiv z^{p+s-1}(c_0 + c_1 z^2+ \dots + c_r z^{2r}).
$$
From this identity we conclude that if $s \neq 2$ then $a_1-b_1
\equiv 0 \mbox{ modulo } p$. If $s=2$, then $a_1-b_1 \equiv c_0
\mbox{ modulo } p$.\\
If $L_{p0}$ has two components then by the Hoste Formula we know
that $c_0$ is minus the  linking number of the two components. As
the link $\overline L_{0}$ has two components (because $\overline
L_{+}$ is a knot), then the linking number of the two components
of $L_{p0}$ is $p$ times the linking number of the two components
of the link $\overline L_{p
0}$. Hence $c_0$ is zero modulo $p$. \\
Consequently, $L_{p+}$ and $L_{p-}$ have the same $a_1$
 coefficient modulo $p$. Notice, that starting from a strongly
 $p-$periodic link $L$, one may change crossings along orbits  to get
 an algebraically split link. The algorithm can be described as
 follows. We start by changing the crossings of $l_1^1$ and
 $l_1^2$, until having their linking number zero. Of course, at the same
 moment we get the same for $l_1^2$ and
 $l_1^3$ ... etc. In the second step we do the same for $l_1^1$ and
 $l_1^3$. Obviously by this procedure  we get an
 algebraically  split link. If $p>3$, then using proposition 4.4 we conclude that $a_1(L)$ is zero
 modulo $p$. It remains to check the case $p=3$. \\
 For $1\leq i \leq 3$, let $\sum_{i}$ be a Seifert surface which is
 bounded by the component $l_1^{i}$. Following the notations in
paragraph 4, recall that  $\bar \mu_{ijk}$ is defined to be  the
algebraic sum of the triple
 points of $\sum_i \cap \sum_j \cap \sum_k $ if $i,j$ and $k$ are
 distinct and  to be zero otherwise.  Now let
 $u_{ij}=\sum_{k=1}^3\bar\mu_{ijk}$. It is well known \cite{Le},  that
 $a_{1}$ is the determinant of the matrix $u_{ij}$.
In our case we can choose  three Seifert surfaces such  that
 $\varphi(\sum_1)=\sum_2$, $\varphi(\sum_2)=\sum_3$ and
 $\varphi(\sum_3)=\sum_1$.
Assume that
 $M \in \sum_i \cap \sum_j \cap \sum_k$, then $\varphi(M)$ and
 $\varphi^2(M)$ are also in the intersection of the three surfaces.
 If $M \notin \Delta$ then the three points are distinct and the
 sum gives zero modulo three.  If $M\in \Delta$, then as the
 linking number of the factor knot with the axis is
 zero modulo three, then the number of intersections points (with signs ) of a
 the quotient surface $\overline {\sum}$ and the axis $\Delta$ is a multiple of 3. It is
 not hard  to see  that the matrix $u_{ij}$, is null modulo 3.
 Hence the coefficient $a_1$ is null modulo 3.\\

Suppose now that  Theorem 2  is true for strongly periodic OS
links with a number of components  less than or equal to
$p\alpha$. Let $L$ be a strongly $p$-periodic OS link having
$p(\alpha+1)$ components. It is easy
to see that we have essentially two cases.\\
%(The other  cases can be ignored).\\
 Case 1: If we change a mixed crossing of the factor link then the resulting link
 $L_{p0}$ will be a strongly $p$-periodic OS link with $p\alpha$
 components. If we write lemma 5.1 in this case we can see that:
 $$ z^{p(\alpha +1)-1}(a_0 + a_1 z^2+ \dots + a_m z^{2m})-z^{{p(\alpha +1)-1}}(b_0 + b_1 z^2+ \dots + b_k
z^{2k})\equiv z^{p+p\alpha-1}(c_0 + c_1 z^2+ \dots + c_r z^{2r}).
$$
Hence, $a_1 - b_1\equiv c_1 \mbox { modulo } p$. By the induction
hypothesis the coefficient $c_1$ is zero modulo $p$.\\

Case 2: If we change a self crossing of the factor link. There are two possible situations:\\
 First, If the resulting link $L_{p0}$ has more components than
 $L_{p+}$, in this case we can see immediately  from lemma 5.1 that $a_1 \equiv   b_1\mbox { modulo }
 p$.\\

 Second,  if the
 resulting link $L_{p0}$ has $p\alpha + 2$ components. In this case only one orbit is affected and   the link $L_{p0}$ is
made up of two invariant components $K_1$ and $K_2$ and a strongly
$p$-periodic link with $p\alpha$ components.
$$L_{p0}=K_1 \cup K_2 \cup l_2^1 \cup\dots \cup l_2^p \cup l_3^1 \cup\dots \cup l_3^p\dots \cup l_{(\alpha+1)}^1 \cup\dots
\cup l_{(\alpha+1)}^p.$$ From lemma 5.1 we have:
$$ z^{p(\alpha +1)-1}(a_0 + a_1 z^2+ \dots + a_m z^{2m})-z^{{p(\alpha +1)-1}}(b_0 + b_1 z^2+ \dots + b_k
z^{2k})\equiv z^{p+p\alpha+1}(c_0 + c_1 z^2+ \dots + c_r z^{2r}).
$$
Consequently, $a_1 - b_1\equiv c_0 \mbox { modulo } p$. It remains
now to prove that the coefficient $c_0$ vanishes modulo $p$.\\
The coefficient $c_0$ can be computed using the Hoste formula
\cite{Ho}. Let us now compute the coefficient $c_0$ in the case of
our link $L_0$. Recall that we assumed that we only changed the
first orbit $l_1^1 \cup\dots \cup l_1^p$ of the link $L$. This
sub-link is transformed onto  the link $K_1\cup K_2$. It can be
easily seen that each of $K_1$ and $K_2$ is invariant by the
rotation. Let $\overline K_1$ and $\overline K_2$ the
corresponding factor knots. We have the following
 $Lk(K_1,K_2)=pLk(\overline
K_1,\overline K_2)\equiv 0 \mbox {
 modulo } p$ and
 $$
Lk(K_i,l_j^1)=Lk(K_i,l_j^2)=\dots =Lk(K_i,l_j^p)\;\; \forall \;\;
1\leq i \leq 2 \mbox{ and } 2\leq j \leq \alpha+1.$$ Furthermore,
as the link $L$ is orbitally separated we have
$\displaystyle\sum_{t=1}^{p} Lk(l_1^i,l_j^t)=0$ for all $i$ and
all $j\neq 1$. Thus, $\displaystyle\sum_{t=1}^{p}
Lk(K_1,l_j^t)+\displaystyle\sum_{t=1}^{p}Lk(K_1,l_j^t)=0$. From
this we conclude that for all $2 \leq t \leq \alpha +1$ we have
$pLk(K_1,l_j^t)+pLk(K_2,l_j^t)$ is zero. Hence,
$Lk(K_1,l_j^t)=-Lk(K_2,l_j^t)$. Consequently the linking matrix of
$L_{p0}$ is of the form (here coefficients are considered modulo
$p$).
$$\left (\begin{array}{cccccc}
0&0& t_2^1&t_2^2&\dots&t_{\alpha+1}^p\\
0&0& -t_2^1&-t_2^2&\dots&-t_{\alpha+1}^p\\
t_2^1&-t_2^2&.&.&\dots&.\\
t_2^2&-t_2^2&.&.&\dots&.\\
.&.&.&.&\dots&.\\
.&.&.&.&\dots&.\\
.&.&.&.&\dots&.\\
t_{\alpha+1}^p&-t_{\alpha+1}^p&.&.&\dots&.

\end{array} \right )
$$

We know that in the linking matrix of a link the sum of lines is
zero. Moreover in our matrix the first line and the second line
are dependant. Thus, all cofactors of the matrix are zero.
Consequently, $L_{p+}$ and $L_{p-}$ have the same $a_1$
 coefficient modulo $p$.\\
 {\bf Remark. } {\sl  As explained above, starting from a strongly
$p$-periodic OS link $L_{p+}$, one may change crossings along an
orbit without modifying the coefficient  $a_1$ modulo $p$. It is
worth mentioning that to prove this result we only used the fact
that the link $L_{p0}$ is OS. Hence, we may apply the same
procedure again
to the link $L_{p-}$, although $L_{p-}$ is not orbitally separated.}\\

In order to transform  the link $L$ onto an algebraically split
link, we apply the same techniques as in the case $\alpha
 =1$.
 Obviously by this procedure all orbits will be made up
 of algebraically  split sub-links. Now, we apply the same techniques for components belonging
 to different orbits until getting an algebraically  split link.
 Finally, by using Proposition 4.4 we conclude that $a_1(L)$ is zero
 modulo $p$. This ends the proof of Theorem 2.\\

%%%%%%%%%%%%%%%%%%%%%%%%%%%%%%%%%%%%%%%%%%%%%%%%%%%%%%%%%%%%%%%%%%%%%%%%%%%%%%%%%%%%%%%%%%%%%%%%%%%%%%%
\begin{center}
{\sc VI- Lescop's Surgery Formula}\\
\end{center}
%%%%%%%%%%%%%%%%%%%%%%%%%%%%%%%%%%%%%%%%%%%%%%%%%%%%%%%%%%%%%%%%%%%%%%%%%%%%%%%%%%%%%%%%%%%%%%%%%%%%%%%%
In \cite{Les}, Lescop provided a combinatorial  description of the
Casson-Walker invariant. Indeed, Lescop introduced  three
equivalent formulas relating the Casson-Walker invariant to the
Conway polynomial and the multi-variable Alexander polynomial.
This allowed her to extend the Casson-Walker invariant to all
oriented three-manifolds. In this paragraph we roughly recall some
notations, we refer the reader to \cite{Les}, ( page 19 and 20)
for more
details.\\
Throughout  the rest of this paper, $N=\{1,\dots,n\}$. Let
$L=\displaystyle\cup_1^n l_i$ an $n$ component link in the
three-sphere. For every $I\subset N$, we define $L_I=\cup_{i\in
I}l_i$. Let $G$ be a graph with $n$ vertices $\{v_1,\dots,v_n\}$
and $e_{ij}$ the edge of $G$ with end points $v_i$ and $v_j$. The
linking number of $L$ with respect to the edge $e_{ij}$,
Lk$(L,e_{ij})$ is defined to be the linking number of $l_i$ and
$l_j$. The linking of $L$ with respect to the graph $G$ is defined
by:
$$\mbox{Lk}(L,G)=\prod_{e \mbox{ \begin{footnotesize}{edge of }\end{footnotesize}G}} \mbox{Lk}(L,e).$$
Let $\cal C$ be the set of graphs $G$ with $n$ vertices such that
$G$ is topologically an oriented circle. We define the circular
linking number of the link $L$ by the following:
$$ \mbox{Lk}_c(L)=\displaystyle\sum_{C\in {\cal C}}\mbox{Lk}(L,C).$$

Now we define the integer  $\theta_b$ as follows:

$$\theta_b(L_I)=\displaystyle\sum_{ K\subset I} \mbox { Lk}_c(L_K) \displaystyle\sum_{(i,j)\in K^2}
\displaystyle\sum_{g\in \sigma_{I \setminus K}}
l_{ig(1)}l_{g(1)g(2)}....l_{g(\sharp (I\setminus K))j}
$$
where $\sigma_{I\setminus K}$ denotes the set of bijections from
$\{1,\dots,\sharp I \setminus K\}$ to $I \setminus K$. Recall here
that we consider surgery with integral coefficients, hence the
formula 1.7.6 in \cite{Les} can be written as follows
$$ \theta(L_I)= \left \{
\begin{array}{ll}
\theta_b(L_I) &\mbox { if } I > 2 \\
\theta_b(L_I)-2l_{ij}& \mbox { if } I =\{i,j\}\\
\theta_b(L_I)+2& \mbox { if } I =\{i\}
\end{array}
\right.
$$
If $I \subset \{1,2,\dots,n\}$, we define $E(L_{I})$ as the
linking matrix of $L_I$ and we define $E(L_{N\setminus
I};I)=(l_{ijI})_{i,j \in N\setminus I}$ where:

$$
l_{ijI}= \left \{
\begin{array}{ll}
l_{ij}& \mbox { if } i \neq j \\
l_{ii}+\displaystyle\sum_{k \in I} l_{ik}& \mbox { if } i=j
\end{array}
\right. $$ Let $M$ be a three manifold obtained from $\S$ by
integral surgery along a framed link $L$ in $\S$. Let $E(L)$ be
the symmetric linking matrix of $L$. Let $b_{-}(L)$ be the number
of negative eigenvalues of $E(L)$ and define $sign(L)$ to be
$(-1)^{b_{-}(L)}$. With respect to this  notations, the
Casson-Walker invariant of $M$ is given by the following formula.

$$\begin{array}{rll}
\lambda(M)=&sign(L)\displaystyle\sum_{J\subset N/ J\neq
\emptyset}det(E(L_{N \setminus J},J))a_1 (L_J) &+\\
&sign(L)\displaystyle\sum_{J\subset N/ J\neq \emptyset}
\frac{det(E(L_{N \setminus J}))(-1)^{\sharp J}\theta (L_J)}{24}&+\\
& |H_1(M)|\displaystyle\frac {Sig(E(L))}{8}.&
\end{array}
$$
This formula corresponds to Proposition 1.7.8 in \cite{Les}. As we
only consider integral surgery then for all $i$ we have $q_i=1$
and the Dedekind sums are null.

%%%%%%%%%%%%%%%%%%%%%%%%%%%%%%%%%%%%%%%%%%%%%%%%%%%%%%%%%%%%%%%%%%%%%%
\begin{center}
{\sc VII- Proof of the main theorem}
\end{center}
%%%%%%%%%%%%%%%%%%%%%%%%%%%%%%%%%%%%%%%%%%%%%%%%%%%%%%%%%%%%%%%%%%%%%%

Let $M$ be a $p$-periodic 3-manifold and $L$ a strongly
$p$-periodic surgery presentation of $L$.
 Let:

$$ \begin{array}{l}
D_1=sign(L)\displaystyle\sum_{J\subset N/ J\neq
\emptyset}det(E(L_{N \setminus J},J))a_1 (L_J)\\
D_2=sign(L)\displaystyle\sum_{J\subset N/ J\neq \emptyset}
\frac{det(E(L_{N \setminus J}))(-1)^{\sharp J}\theta (L_J)}{24}.
\end{array}
$$

{\bf Lemma 7.1.} {\sl Let $p$ be an odd
 prime and $L$ an OS
strongly
$p$-periodic link. Then $D_1$ is null modulo $p$.}\\
 {\bf Proof.} The summation is done over all sub-links of $L$. It
 is easy to see that if the sub-link  $L_J$ is not periodic then there
 is $p$ sub-links of $L$ which are the same as $L_J$. Hence, the
 contributions of these links add to zero modulo $p$. However if
 the sub-link is $p$-periodic, it is necessary an  OS strongly
 $p$-periodic link. By the  Theorem 2, the second
 coefficient of the Conway polynomial of $L_J$ is null modulo $p$. This ends the proof of lemma 7.1.\\

 {\bf Lemma 7.2.} {\sl Let $p$ be an odd prime and $L$ an OS strongly
$p$-periodic link. Then $24D_2$ is null modulo $p$.}\\
{\bf Proof.} In the term $D_2$, the sum is taken over all
sub-links of $L$. As in the proof of the previous lemma, the
contributions of non periodic sub-links add to zero modulo $p$.
Hence, we have to consider only strongly $p$-periodic sub-links.
Recall that if $L_I$ is a link with at least 3 components, then
$$\theta(L_I)=\theta_b(L_I)=\displaystyle\sum_{ K\subset I} \mbox{Lk}_c(L_K)
\displaystyle\sum_{(i,j)\in K^2} \displaystyle\sum_{g\in
\sigma_{I\setminus K}} l_{ig(1)}l_{g(1)g(2)}....l_{g(\sharp (I
\setminus K))j}.
$$
Let $L_J$ a periodic sub-link of $L$. Here also, for the
computation of $\theta(L_J)$ modulo $p$, we need to  consider only
periodic sub-links of $L_J$. If $L_K$ is a periodic sub-link, then
the term $l_{ig(1)}l_{g(1)g(2)}....l_{g(\sharp (I \setminus K))j}$
does not change if the indexes are cyclically permuted. Thus, we
conclude
that $24D_2$ is null modulo $p$. This ends the proof of Lemma 7.2.\\
The proof of Corollary 2 is based on the fact that if $L$ is an
algebraically split $p$-strongly periodic then the linking matrix
of $L$ is diagonal and each diagonal coefficient of the matrix
appears exactly $sp$ times, for some integer $s$. Hence, the
signature of $E(L)$ is null modulo
$p$. The proof of Corollary 1 is obvious.\\
%%%%%%%%%%%%%%%%%%%%%%%%%%%%%%%%%%%%%%%%%%%%%%%%%%%%%%%%%%%%%%%%%%%%%%%%%%%%%%%%%%%%
%%%%%%%%%%%%%%%%%%%%%%%%%%%%%%%%%%%%%%%%%%%%%%%%


\begin{thebibliography}{99}
\bibitem{BGP} J. K. {\sc Bartoszynska}, P. {\sc Gilmer}
and J. {\sc Przytycki}. {\em 3-Manifold invariants and periodicity
of homology spheres}. Algebraic and Geometric Topology, Volume 2
pp. 825-842, (2002).

\bibitem{BM} H. {\sc Bass} and J. W. {\sc Morgan}. {\em The Smith
conjecture}. Pure and App. Math. 112, New York Academic Press
(1994).

\bibitem{Bu} G. {\sc Burde}. {\em \"{U}ber periodische knoten}. Arkiv der
math. (Basel) 30, pp. 487-492, 1978.


\bibitem{Ch1} N. {\sc Chbili}. {\em Les invariants $\theta_{p}$ des
3-vari\'et\'es p\'eriodiques}. Annales de l'institut Fourier,
Fascicule 4, pp.  1135-1150 (2001).

\bibitem{Ch2} N. {\sc Chbili}. {\em Quantum invariants and finite group actions on
 3-manifolds}, Topology
 Appl. vol 136/1-3 pp. 219-231, 2004.

\bibitem{CL} Q. {\sc Chen} and {\sc  T. Le}. {\em Quantum
invariants of periodic links and periodic 3-manifolds}. Preprint.

\bibitem{Da} A. {\sc Davidow}. {\em Casson's invariant  and twisted double
knots}. Topology Appl. 58 pp. 93-101, 1994.

\bibitem{DL} {\sc J. F. Davis, C. Livingston}. {\em
Alexander polynomials of periodic knots}. Topology, 30, pp.
551-564, 1991.

\bibitem{GK} S. {\sc Garoufalidis} and A. {\sc  Kricker}. {\em Finite type invariants of Cyclic branched
covers}, Preprint.

\bibitem{Gi} P.  {\sc Gilmer}. {\em Quantum invariants of periodic three-manifolds}. Geometry and
Topology monographs, Vol. 2: proceeding of the Kirbyfest, pp.
157-175  (1999).

\bibitem{Ha} R. {\sc Hartley}. {\em The Conway potential function for links.} Comment. Math.
Helv. 58 (1983), no. 3, 365--378.

\bibitem{Hi} J. {\sc Hillman}. {\em New proofs of two theorems on periodic
knots}.  Archiv. Math.  37, pp. 457-461, 1981.

\bibitem{Hoso} F. {\sc Hosokawa}. {\em On $\nabla $-polynomials of links.} Osaka
Math. J. 10 1958 273--282. d

\bibitem{Ho} J.{\sc Hoste}. {\em The first
coefficient of the Conway polynomial}. Proc. Amer. Math. Soc., 95,
pp. 299-302, 1985

\bibitem{Ho1} J.{\sc Hoste}. {\em A formula for Casson's invariant}. Trans. Amer. Math. Soc.,297, pp. 547-562,
1986


\bibitem{Is} K. {\sc Ishibe}. {\em The Casson-Walker invariantfor branched cyclic coversof $\S$ branched
over a doubled knot}. Osaka J. Math. 34, pp. 481-495, 1997.

\bibitem{Jo} V. F. R. {\sc Jones}. {\em Hecke algebra
representations of braid groups and link polynomials}. Annals of
Mathematics, Vol. 126, pp. 335-388, 1987.

\bibitem{Ki}  R. {\sc  Kirby}. {\em  The calculus of framed links in $\S$}, Invent. Math. 45 pp. 35-56 (1978).

\bibitem{Les}  C. {\sc Lescop}. {\em Global surgery formula for the Casson-Walker invariant}, Annals of
 Mathematics Studies, Princeton Univ.
 Press (1996).

\bibitem{Le} J. {\sc
Levine}. {\em The Conway polynomial of an algebraically split
link}. Proceeding of knots 96, edited by S. Suzuki, World
Scientific Publishing Co, pp. 23-29, 1997.

\bibitem{Li1}    W. B. R. {\sc Lickorish}. {\em A representation of orientable combinatorial 3-manifolds}.
Ann. Math. 76 (1962), 531-540.

\bibitem{Mi} Y. {\sc Miyazawa}. {\em Conway polynomials of periodic links}. Osaka J. math.  31, pp. 147-163, 1994.
\bibitem{Mil} J.   {\sc Milnor}. {\em Singular points and complex hypersurfaces}. Ann. Math. studies,
Princeton University Press. 1968.
\bibitem{Mul} D. {\sc Mullins}. {\em The generalized Casson invariant for 2-fold branched covers of $\S$
 and the Jones polynomial}. Topology, 32 419-438, 1993.

\bibitem{Mu} K. {\sc Murasugi}. {\em On periodic knots}. Comment.
Math. Helv. 46, pp. 162-174, 1971.


\bibitem{Pr} J. H. {\sc  Przytycki}. {\em On Murasugi's and Traczyk's criteria for periodic
links}. Math. Ann., 283, pp. 465-478, 1989

\bibitem{Pr1} J. H. {\sc  Przytycki}. {\em An elementary proof of
the Traczyk-Yokota criteria for periodic knots.} Proc. Amer. Math.
Soc., 123, pp. 1607-1611, 1995.

\bibitem{PS} J.  {\sc Przytycki }
and M. {\sc Sokolov}. {\em Surgeries on periodic links and
homology of periodic 3-manifolds.} Math. Proc. Cambridge Phil.
Soc. Vol. 131(2),  pp 295--307. 2001.

\bibitem{Ro} D. {\sc Rolfsen}. {\em Knots and Links}.
Mathematics Lecture Series, Publish or Perish, Inc. Houston,
Texas, 1990.

\bibitem{Sa} M. {\sc Sakuma}. {\em On the polynomials of  periodic links}. Math. Ann.  257, pp. 487-494, 1981.

\bibitem{Tr1} P. {\sc Traczyk}. {\em A criterion for knots
of period 3}. Topology appl, 36, pp. 275-281, 1990.

\bibitem{Tr2} P. {\sc Traczyk}. {\em Periodic knots and the
skein polynomial}. Invent. Math, 106(1), pp. 73-84, 1991.

\bibitem{Wa} K. {\sc Walker}. {\em An extension of Casson's
invariant.} Ann. of Math. Studies 126, Princeton University Press,
1992.

\bibitem{Yo} Y.{\sc  Yokota}. {\em The skein polynomial of
periodic knots}. Math. Ann. 291, pp. 281-291, 1991.





\end{thebibliography}
\end{document}